# À propos d'une nouvelle interprétation du
# *Caput variationis* de Leibniz

Godofredo Iommi Amunátegui

Instituto de Física

Pontificia Universidad Católica de Valparaíso

giommi@ucv.cl

Summary

The concept of *caput variationis* plays a fundamental role in Leibniz's *Dissertatio de Arte Combinatoria*, published in 1666. Hereafter we propose an interpretation of the *caput* based on the notion of "induced character" of a group introduced par Georg Frobenius in 1900.

## I

Il y a déjà plus d'un siècle, Georg Frobenius établissait une formule pour calculer les caractères du groupe symétrique $S_n$[1]. Dans le même article, on trouve un algorithme pour obtenir les caractères induits, par les sousgroupes de Young de $S_n$. Ici nous proposons une nouvelle interprétation, dont le travail de Frobenius est le fondement, du concept de *caput variationis*, défini par Leibniz dans sa *Dissertatio de Arte Combinatoria*[2]. Il faut, certes, souligner le fait que notre démarche attribue à Leibniz des développements qui lui sont inconnus et, par ailleurs, fort postérieurs. Toutefois, cette attribution nous semble pertinente car elle peut être entreprise sans faire violence à ses textes. D'emblée il convient de dire que ces pages s'appuient, en partie, sur l'oeuvre de E. Knobloch[3]. En effet quiconque s'intéresse à la combinatoire –au sens large– du philosophe, se doit d'étudier son apport car il explore et expose de manière remarquable cet aspect de l'opus leibnizien.

## II

Dans la *Dissertatio*, Leibniz définit ainsi le *caput:*

*Caput variationis est positio certarum partium*[4]

Le *caput* d'une variation est la position de certaines parties.

Cette définition permet de dire que le *caput* est une caractéristique spatiale. Le versant combinatoire du concept comparaît dans le septième problème, véritable pivot théorique de l'œuvre:

*Dato capite variationes reperire*[5]

Étant donné le *caput* trouver les variations.

Dans ses travaux, Knobloch revient à plusieurs reprises sur le thème du *caput*:

- "Problem 7 considers permutations that contain a *caput*, i.e., a subset that is mapped onto itself by the permutation, or in a special case, remains invariant"[6].
- "A *caput* is a given subset of a set. Leibniz asked for the number of variations that contain such a subset ('*caput variationis*')[7].

Permutations, ensemble et sous ensemble sont les mots-clé de ces énoncés. Or, il s'avère possible d'envisager la *Dissertatio* à partir de la théorie des groupes[8]. Aussi est-ce dans le langage de cette théorie que nous essayerons de répondre à la question posée par Knobloch et issue d'une analyse poussée de la Dissertatio:

> "… Leibniz gives a general solution of the following problem: how to find all possible combinations or all combinations of a given size for a given *caput*. Since he defines a *caput* as a definite subset of definitely given elements which have to be contained in the desired combinations, the problem reads in modern terms: how many combinations of a certain size or of all possible sizes contain a certain number of given elements?"[9]

Désignons cette question par le lettre $Q_1$.

### III

Pour rendre compte du résultat de Frobenius qui nous interesse, il est nécessaire d'exposer quelques détails techniques, en suivant le texte d'aussi près que possible. L'on sait que les permutations de n symboles sont les éléments du groupe symétrique d'ordre n! Les permutations dont la structure cyclique est la même forment une classe. L'ordre de la classe ρ qui a α cycles d'ordre 1, β cycles d'ordre 2, γ cycles d'ordre 3, … est donné par

$$h_\rho = \frac{n!}{1^\alpha \alpha! \, 2^\beta \beta! \, 3^\gamma \gamma! \, \ldots}$$

Le nombre k des classes est le nombre de solutions entières positives de l'équation n = α + 2β + 3γ+ … Le nombre k des classes est égal au nombre des caractères du groupe symétrique $S_n$. L'on a ρ = 0, 1, … k-1. Par ailleurs n peut s'écrire ainsi: n = $n_1$ + $n_2$ + $n_3$ + …Cette somme correspond à une partition de n. Les sous-partitions $n_1$, $n_2$, … sont des entiers positifs ou zéro. Chaque partition de n constitue, à son tour, un sousgroupe de $S_n$ dont l'ordre es g = $n_1$!$n_2$! … Chaque permutation du sousgroupe $S_{n_1} \cdot S_{n_2} \cdot \ldots$ se décompose en une permutation de $n_1$ symboles, en une permutation des $n_2$ symboles etc. Si l'on considère la permutation $R_\upsilon$ de $\alpha_\upsilon$ cycles de longueur 1, $\beta_\upsilon$ cycles de longueur 2 etc. Il vient: $n_1$ = $\alpha_1$ + 2$\beta_1$ + …, $n_2$ = $\alpha_2$ + 2$\beta_2$ + … D'autre part α = $\alpha_1$ + $\alpha_2$ + …, β = $\beta_1$ + $\beta_2$ + …

Le nombre d'éléments de $S_{n_1}, S_{n_2}, \ldots$ est le produit des ordres des classes de $S_{n_1}, S_{n_2}$, etc. C'est-à-dire:

$$\frac{n_1!}{1^{\alpha_1}\alpha_1! 2^{\beta_1}\beta_1!} \cdots \frac{n_2!}{1^{\alpha_2}\alpha_2! 2^{\beta_2}\beta_2!} \ldots$$

Le nombre total des éléments de la classe ρ en $S_{n_1}, S_{n_2}, \ldots$ est donné par la somme

$$g_\rho = \sum \frac{n_1!}{1^{\alpha_1}\alpha_1! 2^{\beta_1}\beta_1!} \cdots \frac{n_2!}{1^{\alpha_2}\alpha_2! 2^{\beta_2}\beta_2!} \cdots$$

qui comprend toutes les partitions.

Frobenius construit l'expression générale pour le caractère induit:

$$\frac{g}{h} \frac{h_\rho}{g_\rho} = \sum \frac{\alpha!}{\alpha_1! \alpha_2!} \cdots \frac{\beta!}{\beta_1! \beta_2!} \cdots$$

Dans laquelle: g = ordre du groupe G

$g_\rho$ = ordre de la classe ρ de G.

h = ordre du sousgroupe H de G.

$h_\rho$ = ordre de la classe ρ de H, i.e. cette classe comprend $h_\rho$ éléments de la classe ρ.

Dans la notation actuelle (λ) = $(\lambda_1, \ldots \lambda_p)$ est une partition de n où n = $\lambda_1 + \lambda_2 + \ldots \lambda_p$, $\lambda_1 \geq \lambda_2 \geq \ldots \geq \lambda_p \geq 0$; p(n) est le nombre des partitions de n. Chaque partition de n est liée au sousgroupe $S_\lambda$ de $S_n$, donné par le produit direct $S_{\lambda_1} \times S_{\lambda_2} \times S_{\lambda_3} \times \cdots \times S_{\lambda_p}$. ces $S_{\lambda_i (i=1,-p)}$ sont les sousgroupes de Young de $S_n$ (Frobenius écrit $S_{n_1}, S_{n_2}, \ldots$ ).

Exemple: calcul des caractères induits par le sousgroupe $S_3 \times S_2$ de $S_5$.

(i) g: ordre de $S_5$ = 5! = 120.
(ii) Ordre des classes de $S_5$: (5) = 24; (4,1) = 30; (3,2) = 20; (3,$1^2$) = 20; ($2^2$,1) = 15, (2,$1^3$) = 10, ($1^5$) = 1.

(iii)    h: ordre de $S_3 \times S_2$ = 3!2! = 12.

(iv) Les éléments de $S_3$ et de $S_2$, écrits en notation cyclique, sont respectivement: {(·)(·)(·); 3(·)(··); 2(···)} et {(·)(·); (··)}, donc les éléments de $S_3 \times S_2$ sont: {(·)(·)(·)(·)(·); 4(··)(·)(·)(·); 2(···)(·)(·); 2(···)(··); 3(·)(··)(··)}. Notons que ces éléments ne comportent nul cycle d'ordre 4 ou 5, en conséquence les caractères induits correspondant aux partitions (4,1) et (5) sont 0. On notera $\phi^{\lambda}_{(1^{\alpha},2^{\beta}...)}$ le caractère induit correspondant au sousgroupe de Young $S_{\lambda_1} \times S_{\lambda_2} \times \cdots \times S_{\lambda_p}$ et à la clase $(1^{\alpha}, 2^{\beta}, 3^{\gamma} \ldots)$

(v) Donc, il vient:

$$\phi^{3,2}_{1^5} = \frac{120}{12} \cdot \frac{1}{1} = 10; \quad \phi^{3,2}_{2,1^3} = \frac{120}{12} \cdot \frac{4}{10} = 4; \quad \phi^{3,2}_{2^2,1} = \frac{120}{12} \cdot \frac{3}{15} = 2$$

$$\phi^{3,2}_{3,1^2} = \frac{120}{12} \cdot \frac{2}{20} = 1 \quad et \quad \phi^{3,2}_{3,2} = \frac{120}{12} \cdot \frac{2}{20} = 1$$

En bref:

*Caractères induits par $S_3 \times S_2$ de $S_5$*

|       | $(1^5)$ | $(2,1^3)$ | $(2^2,1)$ | $(3,1^2)$ | $(3,2)$ | $(4,1)$ | $(5)$ |
|-------|---------|-----------|-----------|-----------|---------|---------|-------|
| (3,2) | 10      | 4         | 2         | 1         | 1       | 0       | 0     |

En général, les caractères forment une matrice carrée de dimensions p(n) × p(n) dont les files et les colonnes sont indexées par les partitions de n arrangées selon l'ordre lexicographique, et par les classes. A titre d'illustration, voici la matrice pour n = 5.

*Classes*

|     | $1^5$ | $1^3,2$ | $1,2^2$ | $1^2,3$ | $2,3$ | $1,4$ | $5$ |
|-----|-------|---------|---------|---------|-------|-------|-----|
| 5   | 1     | 1       | 1       | 1       | 1     | 1     | 1   |
| 4,1 | 5     | 3       | 1       | 2       | 0     | 1     | 0   |
| 3,2 | 10    | 4       | 2       | 1       | 1     | 0     | 0   |
| $3,1^2$ | 20 | 6 | 0 | 2 | 0 | 0 | 0 |
| $2^2,1$ | 30 | 6 | 2 | 0 | 0 | 0 | 0 |
| $2,1^3$ | 60 | 6 | 0 | 0 | 0 | 0 | 0 |
| $1^5$ | 120 | 0 | 0 | 0 | 0 | 0 | 0 |

## IV

Reprenons la question $Q_1$. En un premier temps nous allons l'écrire en y introduisant la notation cyclique (ci-après $Q_2$)

"Combien de combinaisons, de toutes les structures cycliques possibles, contiennent un certain nombre d'éléments?"

Décomposons $Q_2$:

- Le "nombre d'éléments" indique les combinaisons qui appartiennent au groupe symétrique $S_n$, déterminées par les partitions de n et par les classes.
- Les "combinaisons de toutes les structures cycliques possibles" peuvent être attribuées aux caractères induits en $S_n$ par les sousgroupes $S_{\lambda_1} \times \cdots \times S_{\lambda_p}$.

Maintenant $Q_2$ peut se réécrire de la façon suivante ($Q_3$):

"Combien d'éléments du groupe symétrique $S_n$, correspondant à la partition de n = $\lambda_1 + \ldots \lambda_p$, et à la clase $(1^\alpha, 2^\beta, 3^\gamma, \ldots)$, sont contenus dans les caractères induits par les sousgroupes $S_{\lambda_1} \times \cdots \times S_{\lambda_p}$ en $S_n$?"

Le nombre cherché est égal à $\phi^\lambda_{(1^\alpha, 2^\beta, 3^\gamma, \ldots)}$.

Partant, le passage de $Q_1$ à $Q_3$ constitue notre interprétation du *caput variationis.*